\documentclass{amsart}
\usepackage{amssymb}
\usepackage{upref}

\usepackage[all,cmtip]{xy}

\usepackage[lite,initials,shortalphabetic]{amsrefs}

\usepackage[colorlinks]{hyperref}

\providecommand{\texorpdfstring}[2]{#1}

\makeatletter
  \renewcommand{\p@enumi}{\thesubsection}
\makeatother

\newenvironment{resumeenumerate}[1]
{\begin{enumerate}
 \setcounter{enumi}{#1}
 \addtocounter{enumi}{-1}
}
{\end{enumerate}
}

\newenvironment{lettered}
{\begin{list}{\thelettercounter.}
 {\usecounter{lettercounter}\def\makelabel##1{\hss\llap{##1}}}
}
{\end{list}
}
\newcounter{lettercounter}
\renewcommand{\thelettercounter}{\Alph{lettercounter}}

\theoremstyle{plain}

\makeatletter
\newcommand{\emsection}[1]{%
  \par
  \addpenalty\@secpenalty
  \vskip 6 pt plus 9 pt
  \emph{#1.}\nobreak\enspace\ignorespaces
}
\makeatother

\newcommand{\intro}{%
  \goodbreak
  \vskip 6 pt plus 9 pt
}

\numberwithin{equation}{subsection}

\newcommand{\Comma}{\rlap{\enspace ,}}
\newcommand{\Period}{\rlap{\enspace .}}

\newcommand{\cat}[1]{\boldsymbol{#1}}

\newcommand{\Cat}{\mathbf{Cat}}

\DeclareMathOperator{\Gr}{\boldsymbol{Gr}}

\newcommand{\op}{^{\mathrm{op}}}

\newcommand{\spacedcdots}{{\cdot\;\cdot\;\cdot}}

\newcommand{\zigzag}[5]{#1\negmedspace: #2
  \mathchoice{\longrightarrow}{\to}{\to}{\to}#3
  \mathchoice{\longleftarrow}{\gets}{\gets}{\gets}#4:\negmedspace #5}

\begin{document}

\title{A Quillen Theorem $B_{n}$ for homotopy pullbacks}

\author{C. Barwick}
\address{Department of Mathematics, Massachusetts Institute of
  Technology, Cambridge, MA 02139}
\email{clarkbar@math.mit.edu}

\author{D.M. Kan}
\address{Department of Mathematics, Massachusetts Institute of
  Technology, Cambridge, MA 02139}

\date{\today}

\begin{abstract}
  We prove an extension of the Quillen Theorem $B_{n}$ for homotopy
  \emph{fibres} of \cite{DKS}*{\S 6} to a similar result for homotopy
  \emph{pullbacks} and use this to obtain sufficient conditions on a
  zigzag $\cat X \to \cat Y \gets \cat Z$ between categories in order
  that \emph{its pullback is a homotopy pullback}.
\end{abstract}

\maketitle

\section{Introduction}
\label{sec:Intro}

\subsection{The background}
\label{sec:Bkgrnd}

In \cite{Q}*{\S 1} Quillen proved his Theorem B which gave a rather
simple description of the homotopy fibres of a functor $f\colon \cat X
\to \cat Y$ if $f$ had a certain \emph{property $B_{1}$}.

This was generalized in \cite{DKS}*{\S 6} where it was shown that
increasingly weaker \emph{properties} $B_{n}$ ($n>1$) allowed for
increasingly less simple descriptions of these homotopy fibres.
Moreover it was noted that a sufficient condition for a functor
$f\colon \cat X \to \cat Y$ to have property $B_{n}$ ($n>1$) was that
the category $\cat Y$ has a certain \emph{property $C_{n}$}.

\subsection{The current paper}
\label{sec:CurPpr}

We show that for a zigzag $\zigzag{f}{\cat X}{\cat Y}{\cat Z}{g}$ in
which $f$ has property $B_{n}$ \eqref{sec:Bkgrnd} (and in particular
if $\cat Y$ has property $C_{n}$ \eqref{sec:Bkgrnd}), its homotopy
pullback admits a description rather similar to the ones mentioned in
\ref{sec:Bkgrnd}.

Moreover the pullback $\cat X\times_{\cat Y}\cat Z$ of this zigzag
comes with a monomorphism into this homotopy pullback and hence is
itself a homotopy pullback if the monomorphism is a weak equivalence.

\subsection{The motivation}
\label{sec:Mtvtn}

Our result \eqref{sec:CurPpr} and in particular its second half is
what really motivated us to write the present note, and well for the
following reasons.

In \cite{R}*{8.3} Charles Rezk proved that
\begin{itemize}
\item for every \emph{simplicial model category} one (and hence every)
  Reedy fibrant replacement of its simplicial nerve is a
  \emph{complete Segal space}.
\end{itemize}

Although the proof of the Segal part of this result relied heavily on
the simplicial structure, it seemed that this result would also hold
without the assumption of a simplicial structure.

In fact as we will show in \cite{BK}, most of the model structure is
superfluous.  All that is needed is that there is a category of weak
equivalences with three rather simple properties.  More precisely we
will show that
\begin{itemize}
\item Charles Rezk's result holds for every \emph{relative category}
  which has the \emph{two out of six property} and admits a
  \emph{$3$-arrow calculus}.
\end{itemize}

It turns out that in that situation the category of the weak
equivalences has property $C_{3}$ with the result that, in view of our
result \eqref{sec:CurPpr} for $n = 3$, the verification of the Segal
property, i.e.\ showing that certain fibre products (which are
iterated pullbacks) are homotopy fibre products (which are iterated
homotopy pullbacks), is reduced to a rather simple calculation.

\subsection{The proof}
\label{sec:Prf}

The homotopy fibre results of \eqref{sec:Bkgrnd} were obtained by an
induction on $n$ which at each stage used Quillen's Theorem~B.

To prove our homotopy pullback results \eqref{sec:CurPpr} it turns out
to be convenient to go one step further back to the lemma that Quillen
used to prove his Theorem~B and which can be summarized as follows:
\begin{itemize}
\item If $F\colon \cat D \to \Cat$ is a $\cat D$-diagram of categories
  and \emph{weak equivalences} between them, $\Gr F$ its Grothendieck
  construction and $\pi\colon \Gr F \to \cat D$ the associated
  projection functor, then, \emph{for every object $D \in \cat D$, the
    fibre}
  \begin{displaymath}
    \pi^{-1}D = F D
  \end{displaymath}
  \emph{of $\pi$ over $D$ is also a homotopy fibre}.
\end{itemize}

Using this result we first give a different non-inductive proof of the
results of \eqref{sec:Bkgrnd} and then note that this proof almost
effortless extends to a proof of the homotopy pullback results of
\eqref{sec:CurPpr}.

\subsection{Organization of the paper}
\label{sec:OrgPpr}

There are three more sections.

In the first (\S \ref{sec:prelim}) we discuss various
\emph{Grothendieck constructions} and give a precise formulation of
what we will call \emph{Quillen's lemma}.

In the next section (\S \ref{sec:Rslts}) we recall the properties
$B_{n}$ and $C_{n}$ and state the \emph{Theorems $B_{n}$ for homotopy
  fibres} and for \emph{homotopy pullbacks}.

The last section (\S \ref{sec:Prfs}) then is devoted to a proof of
these two Theorems $B_{n}$.

\section{Preliminaries}
\label{sec:prelim}

In preparation for the formulation and the proofs of our results we
here
\begin{itemize}
\item briefly discuss Grothendieck constructions,
\item formulate, in terms of Grothendieck constructions, a categorical
  version of the lemma that Quillen used in his proof of Theorem B, and
\item describe three Grothendieck constructions which will be used in
  our proofs.
\end{itemize}

But first a comment on
\addtocounter{subsection}{-1}
\subsection{Terminology}
\label{sec:Term}

We will work in the category $\Cat$ of small categories with the
Thomason model structure \cite{T2} in which a map is a \emph{weak
  equivalence} iff its \emph{nerve} is a weak equivalence of
simplicial sets and in which \emph{homotopy fibres} and \emph{homotopy
pullbacks} have a similar meaning.

\subsection{Grothendieck constructions}
\label{sec:GrthCnst}

Given a small category $\cat D$ and a functor $F\colon \cat D \to
\Cat$ \eqref{sec:Term}, the \textbf{Grothendieck construction} on $F$
is the category $\Gr F$ which has
\begin{enumerate}
\item \label{GrthCnst:i} as \emph{objects} the pairs $(D,A)$
  consisting of objects
  \begin{displaymath}
    D \in \cat D
    \qquad\text{and}\qquad
    A \in F D \Comma
  \end{displaymath}
\item \label{GrthCnst:ii} as \emph{maps} $(D_{1},A_{1}) \to (D_{2},
  A_{2})$ the pairs $(d,a)$ of maps
  \begin{displaymath}
    d\colon D_{1}\to D_{2} \in \cat D
    \qquad\text{and}\qquad
    a\colon(F d)A_{1} \to A_{2} \in F D_{2}
  \end{displaymath}
\end{enumerate}
and
\begin{resumeenumerate}{3}
\item \label{GrthCnst:iii} in which the \emph{composition} is given by
  the formula
  \begin{displaymath}
    (d',a')(d,a) = \bigl(d'd, a'\bigl((F d)a\bigr)\bigr) \Period
  \end{displaymath}
\end{resumeenumerate}
Moreover
\begin{resumeenumerate}{4}
\item \label{GrthCnst:iv} $\Gr F$ comes with a \emph{projection
    functor} $\pi\colon \Gr F \to \cat D$ which sends an object
  $(D,A)$ (resp.\ a map $(d,a)$) in $\Gr F$ to the object $D$ (resp.\
  the map $d$) in $\cat D$.
\end{resumeenumerate}

The usefulness of Grothendieck constructions is due to the following
property which was noticed by Bob Thomason \cite{T1}*{1.2}:
\subsection{Proposition}
\label{prop:GrHoCo}

\begin{em}
  \begin{enumerate}
  \item \label{GrHoCo:i} The Grothendieck construction is a homotopy
    colimit construction on the category $\Cat$,
  \end{enumerate}
  and hence
  \begin{resumeenumerate}{2}
  \item \label{GrHoCo:ii} it is homotopy invariant in the sense that
    every natural weak equivalence \eqref{sec:Term} between two
    functors $F_{1},F_{2}\colon \cat D \to \Cat$ induces a weak
    equivalence $\Gr F_{1} \to \Gr F_{2}$.
  \end{resumeenumerate}
\end{em}

\intro
Next we note that Quillen's key observation in the lemma that he used
to prove Theorem B was that certain functors $\cat D \to \Cat$ had
what we will call
\subsection{Property \texorpdfstring{$Q$}{Q}}
\label{sec:PropQ}

Given a small category $\cat D$, a functor $F\colon \cat D \to \Cat$
will be said to have \textbf{property $Q$} if it sends all maps of $\cat
D$ to weak equivalences in $\Cat$.

A categorical version of the lemma that Quillen used in the proof of
Theorem B (a proof of which can be found in \cite{GJ}*{IV, 5.7}) then
becomes in view of \ref{GrHoCo:i} above
\subsection{Quillen's lemma}
\label{sec:QuilLem}

\begin{em}
  If, given a small category $\cat D$, a functor $F\colon \cat D \to
  \Cat$ has property $Q$ \eqref{sec:PropQ}, then, for every object $D
  \in \cat D$, the fibre
  \begin{displaymath}
    \pi^{-1}D = FD
  \end{displaymath}
  of $\pi$ \eqref{GrthCnst:iv} over $\cat D$ is a homotopy fibre.
\end{em}

\intro
It remains to construct the promised three Grothendieck
constructions.

We start with
\subsection{Two Grothendieck constructions associated with a functor
  \texorpdfstring{$\cat X \to \cat Y$}{}}
\label{sec:TwoGroth}

Given an integer $n \ge 1$ and a functor $f\colon \cat X \to \cat Y$
between small categories, we denote by $(f\cat X \downarrow_{n}\cat
Y)$ the category of which
\begin{enumerate}
\item \label{TwoGroth:i} an object consists of a pair of objects
  \begin{displaymath}
    X \in \cat X
    \qquad\text{and}\qquad
    Y \in \cat Y
  \end{displaymath}
  together with an alternating zigzag
  \begin{displaymath}
    fX = Y_{n} \quad\spacedcdots\quad Y_{2}
    \longleftarrow Y_{1} \longrightarrow gZ
    \qquad\text{in $\cat Y$}
  \end{displaymath}
\end{enumerate}
and of which
\begin{resumeenumerate}{2}
\item \label{TwoGroth:ii} a \emph{map} consists of a pair of maps
  \begin{displaymath}
    x\colon X \to X' \in \cat X
    \qquad\text{and}\qquad
    y\colon Y \to Y' \in \cat Y
  \end{displaymath}
  together with a commutative diagram
  \begin{displaymath}
    \vcenter{
      \xymatrix{
        {fX = Y_{n}} \ar@{}[r]|-{\spacedcdots} \ar[d]_{fx}
        & {Y_{2}} \ar[d]
        & {Y_{1}} \ar[l] \ar[r] \ar[d]
        & {Y} \ar[d]^{y}\\
        {fX' = Y'_{n}} \ar@{}[r]|-{\spacedcdots}
        & {Y'_{2}}
        & {Y'_{1}} \ar[l] \ar[r]
        & {Y'}
      }
    }
    \qquad\text{in $\cat Y$}
  \end{displaymath}
\item \label{TwoGroth:iii} This category comes with a monomorphism
  \begin{displaymath}
    h\colon X \longrightarrow (f\cat X \downarrow_{n}\cat Y)
  \end{displaymath}
  which sends each object $X \in \cat X$ to the zigzag of identity
  maps which starts at $fX$.
\end{resumeenumerate}
Furthermore
\begin{resumeenumerate}{4}
\item \label{TwoGroth:iv} let, for every object $Y \in \cat Y$
  \begin{displaymath}
    (f\cat X \downarrow_{n} Y) \subset
    (f\cat X \downarrow_{n}\cat Y)
  \end{displaymath}
  denote the subcategory consisting of the objects which end at $Y$
  and the maps which end at $!_{Y}$
\end{resumeenumerate}
and similarly
\begin{resumeenumerate}{5}
\item \label{TwoGroth:v} let, for every object $X \in \cat X$
  \begin{displaymath}
    (fX\downarrow_{n}\cat Y) \subset
    (f\cat X \downarrow_{n} \cat Y)
  \end{displaymath}
  denote the subcategory consisting of the objects which start at $fX$
  and the maps which start at $1_{fX}$.
\end{resumeenumerate}
The naturality of $(f\cat X \downarrow_{n} Y)$ and
$(fX\downarrow_{n}\cat Y)$ in respectively $Y$ and $X$ then readily
implies
\subsection{Proposition}
\label{prop:oddeven}
\begin{em}
  For every integer $n \ge 1$ and functor $f\colon \cat X \to \cat Y$
  between small categories \eqref{sec:GrthCnst}
  \begin{enumerate}
  \item \label{oddeven:i} $(f\cat X \downarrow_{n}\cat Y) =
    \Gr\bigl((f\cat X \downarrow_{n} -)\colon \cat Y \to \Cat\bigr)$
  \end{enumerate}
  \medskip
  and
  \begin{resumeenumerate}{2}
  \item \label{oddeven:ii} $(f\cat X \downarrow_{n}\cat Y) =
    \left\{
      \begin{gathered}
        \Gr\bigl((f-\downarrow_{n}\cat Y)\colon \cat X \to
          \Cat\bigr)\\
        \text{or}\\
        \Gr\bigl((f-\downarrow_{n}\cat Y)\colon \cat X\op \to
          \Cat\bigr)
      \end{gathered}
    \right.$\\
    depending on whether $n$ is even or odd.
  \end{resumeenumerate}
\end{em}

\intro
We end with
\subsection{A Grothendieck construction associated with a zigzag
  \texorpdfstring{$\cat X \to \cat Y \gets \cat Z$}{}}
\label{sec:GrthZgZg}

Given an integer $n \ge 1$ and a zigzag $\zigzag{f}{\cat X}{\cat
  Y}{\cat Z}{g}$ between small categories, we denote by $(f\cat X
\downarrow_{n} g\cat Z)$ the category of which
\begin{enumerate}
\item \label{GrthZgZg:i} an \emph{object} consists of a pair of objects
  \begin{displaymath}
    X \in \cat X \qquad\text{and}\qquad Z \in \cat Z \Comma
  \end{displaymath}
  together with an alternating zigzag
  \begin{displaymath}
    fX = Y_{n} \quad\spacedcdots\quad Y_{2}
    \longleftarrow Y_{1} \longrightarrow gZ
    \qquad\text{in $\cat Y$}
  \end{displaymath}
\end{enumerate}
and of which
\begin{resumeenumerate}{2}
\item \label{GrthZgZg:ii} a map consists of a pair of maps
  \begin{displaymath}
    x\colon X \to X' \in \cat X
    \qquad\text{and}\qquad
    z\colon Z \to Z' \in \cat Z \Comma
  \end{displaymath}
  together with a commutative diagram
  \begin{displaymath}
    \vcenter{
      \xymatrix{
        {fX = Y_{n}} \ar@{}[r]|-{\spacedcdots} \ar[d]_{fx}
        & {Y_{2}} \ar[d]
        & {Y_{1}} \ar[l] \ar[r] \ar[d]
        & {gZ} \ar[d]^{gz}\\
        {fX' = Y'_{n}} \ar@{}[r]|-{\spacedcdots}
        & {Y'_{2}}
        & {Y'_{1}} \ar[l] \ar[r]
        & {gZ'}
      }
    }
    \qquad\text{in $\cat Y$}\Period
  \end{displaymath}
\item \label{GrthZgZg:iii} This category comes with a monomorphism
  \begin{displaymath}
    K\colon (\cat X \times_{\cat Y} \cat Z) \longrightarrow
    (f\cat X \downarrow_{n} g\cat Z)
  \end{displaymath}
  which sends each object $(X,Z) \in \cat X\times_{\cat Y}\cat Z$ to a
  zigzag of identity maps starting at $fX$ and ending at $gZ$.
\end{resumeenumerate}
Furthermore
\begin{resumeenumerate}{4}
\item \label{GrthZgZg:iv} we denote, for every object $Z \in \cat Z$,
  by
  \begin{displaymath}
    (f\cat X \downarrow_{n}gZ) \subset
    (f\cat X \downarrow_{n}g\cat Z)
  \end{displaymath}
  the subcategory consisting of the objects which end at $gZ$ and the
  maps which end at $1_{gZ}$.
\end{resumeenumerate}
The naturality of $(f\cat X \downarrow_{n}gZ)$ with respect to $Z$
then readily implies
\subsection{Proposition}
\label{prop:ZCztMp}

For every integer $n \ge 1$ and zigzag $\zigzag{f}{\cat X}{\cat
  Y}{\cat Z}{g}$ between small categories \eqref{sec:GrthCnst}
\begin{displaymath}
  (f\cat X\downarrow_{n}g\cat Z) =
  \Gr\bigl((f\cat X\downarrow_{n} g-)\colon \cat Z \to \Cat\bigr)
  \Period
\end{displaymath}

\section{The results}
\label{sec:Rslts}

We start with recalling the homotopy fibre results of
\eqref{sec:Bkgrnd}, beginning with the notion of
\subsection{Property \texorpdfstring{$B_n$}{Bn}}
\label{sec:PropBn}

Given an integer $n \ge 1$, a functor $f\colon \cat X \to \cat Y$
between small categories is said to have \textbf{property $B_{n}$} if
the functor \eqref{oddeven:i}
\begin{displaymath}
  (f\cat X \downarrow_{n} -)\colon \cat Y \longrightarrow \Cat
\end{displaymath}
has property $Q$ \eqref{sec:PropQ}.

\intro
One then has
\subsection{Theorem \texorpdfstring{$B_n$}{Bn}}
\label{thm:Bn}
\begin{em}
  If a functor $f\colon \cat X \to \cat Y$ between small categories
  has property $B_{n}$ ($n\ge 1$), then, for every object $Y \in \cat
  Y$, the category $(f\cat X\downarrow_{n}Y)$ \eqref{TwoGroth:iv} is a
  homotopy fibre of $f$ over $Y$.
\end{em}

\intro
Closely related to property $B_{n}$ is
\subsection{Property \texorpdfstring{$C_n$}{Cn}}
\label{sec:PropCn}

Let $\cat O$ denote the category consisting of a single object and its
identity map and let $n$ be an integer $\ge 1$.  Then a small category
$\cat Y$ is said to have \textbf{property $C_{n}$} if
\begin{enumerate}
\item \label{PropCn:i} every functor $e\colon \cat O \to \cat Y$ has
  property $B_{n}$, i.e.
\item \label{PropCn:ii} every functor $e\colon \cat O \to \cat Y$
  gives rise to a functor $(e\cat O\downarrow_{n}-)\colon\cat Y \to
  \Cat$ which has property $Q$ \eqref{sec:PropQ}.
\end{enumerate}

\intro
The usefulness of this notion is due to the fact that, in view of
\ref{GrHoCo:ii}, \ref{oddeven:ii} and \ref{PropCn:ii}, one has
\subsection{Theorem \texorpdfstring{$C_n$}{Cn}}
\label{thm:Cn}

\begin{em}
  If $f\colon \cat X \to \cat Y$ is a functor between small categories
  and $\cat Y$ has property $C_{n}$ ($n \ge 1$), then $f$ has property
  $B_{n}$.
\end{em}

\intro
We end with formulating our ultimate aim, namely
\subsection{Theorem \texorpdfstring{$B_n$}{Bn} for homotopy pullbacks}
\label{thm:BnPlbk}

\begin{em}
  Let $n$ be an integer $\ge 1$ and let $\zigzag{f}{\cat X}{\cat
    Y}{\cat Z}{g}$ be a zigzag between small categories.  If $f$ has
  property $B_{n}$ \eqref{sec:PropBn} (and in particular if $\cat Y$
  has property $C_{n}$ \eqref{sec:PropCn}), then
  \begin{enumerate}
  \item \label{BnPlbk:i} the category $(f\cat X\downarrow_{n}g\cat Z)$
    \eqref{sec:GrthZgZg} is a homotopy pullback of this zigzag.
  \end{enumerate}
  Moreover if in addition the monomorphism \eqref{GrthZgZg:iii}
  \begin{displaymath}
    k\colon (\cat X\times_{\cat Y}\cat Z) \longrightarrow
    (f\cat X\downarrow_{n}g\cat Z)
  \end{displaymath}
  is a weak equivalence, then
  \begin{resumeenumerate}{2}
  \item \label{BnPlbk:ii} the pullback $(\cat X\times_{\cat Y}\cat Z)$
    of this zigzag is also a homotopy pullback.
  \end{resumeenumerate}
\end{em}

\section{The proofs}
\label{sec:Prfs}

It remains to give a proof of theorems \ref{thm:Bn} and
\ref{thm:BnPlbk}, starting with
\subsection{A proof of Theorem \texorpdfstring{$B_n$}{Bn}
  \texorpdfstring{\eqref{thm:Bn}}{}}
\label{sec:PrfBn}

Given an object $Y \in \cat Y$, it follows from \ref{sec:QuilLem} and
\ref{oddeven:i} that
\begin{enumerate}
\item \label{PrfBn:i} $(f\cat X\downarrow_{n}Y)$ is the fibre as well
  as a homotopy fibre over $Y$ of the projection functor
  \begin{displaymath}
    \pi\colon \Gr(f\cat X\downarrow_{n}-) =
    (f\cat X\downarrow_{n}\cat Y) \longrightarrow \cat Y \Period
  \end{displaymath}
\end{enumerate}
That it is also a homotopy fibre over $Y$ of the functor $f\colon \cat
X \to \cat Y$ therefore is a consequence of
\begin{resumeenumerate}{2}
\item \label{PrfBn:ii} the commutativity of the diagram
  \begin{displaymath}
    \xymatrix{
      {\cat X} \ar[rr]^-{h} \ar[dr]_{f}
      && {(f\cat X\downarrow_{n}\cat Y)} \ar[dl]^{\pi}\\
      & {\cat Y}
    }
  \end{displaymath}
  in which $h$ is as in \ref{TwoGroth:iii}
\end{resumeenumerate}
and the readily verifiable fact that
\begin{resumeenumerate}{3}
\item \label{PrfBn:iii} $h$ is a weak equivalence.
\end{resumeenumerate}

\intro
Finally we are ready to give
\subsection{A proof of Theorem \texorpdfstring{$B_n$}{Bn} for homotopy
  pullbacks \texorpdfstring{\eqref{thm:BnPlbk}}{}}
\label{sec:PrfBnHPB}

As the functor $f\colon \cat X \to \cat Y$ has property $B_{n}$, i.e.
\begin{itemize}
\item the functor $(f\cat X\downarrow_{n}-)\colon\cat Y \to \Cat$ has
  property $Q$ \eqref{sec:PropQ}
\end{itemize}
it readily follows that
\begin{itemize}
\item the functor $(f\cat X\downarrow_{n}g-)\colon \cat Z \to \Cat$
  \eqref{prop:ZCztMp} also has property $Q$.
\end{itemize}
Consequently (\ref{sec:QuilLem} and \ref{prop:ZCztMp})
\begin{enumerate}
\item for every object $Z \in \cat Z$, $(f\cat X\downarrow_{n}gZ)$ is
  the fibre as well as the homotopy fibre over $Z$ of the projection
  functor
  \begin{displaymath}
    \Gr\bigl(f\cat X\downarrow_{n}g-\bigr) =
    (f\cat X\downarrow_{n}g\cat Z) \longrightarrow \cat Z
  \end{displaymath}
  Now consider the commutative square
  \begin{displaymath}
    \xymatrix{
      {(f\cat X\downarrow_{n}\cat Y)} \ar[d]_{\pi}
      & {(f\cat X\downarrow_{n}g\cat Z)} \ar[d]^{\pi} \ar[l]_{g'}\\
      {\cat Y}
      & {\cat Z} \ar[l]_{g}
    }
  \end{displaymath}
  in which $g'$ is induced by $g$.
\end{enumerate}
Then clearly
\begin{resumeenumerate}{2}
\item this square is a pullback square and hence, for every object $Z
  \in \cat Z$, $g'$ maps the fibre over $Z$ isomorphically onto the
  fibre over $gZ \in \cat Y$.
\end{resumeenumerate}
Therefore, in view of (i) and \ref{PrfBn:i}
\begin{resumeenumerate}{3}
\item this pullback square is a homotopy pullback square.
\end{resumeenumerate}
With other words $(f\cat X\downarrow_{n}g\cat Z)$ is a homotopy
pullback of the zigzag
\begin{displaymath}
  \zigzag{\pi}{(f\cat X\downarrow_{n}\cat Y)}{\cat Y}{\cat Z}{g}
\end{displaymath}
That it is also a homotopy pullback of the zigzag
\begin{displaymath}
  \zigzag{f}{\cat X}{\cat Y}{\cat Z}{g}
\end{displaymath}
now follows from \ref{PrfBn:ii} and \ref{PrfBn:iii}.


\end{document}